\documentstyle[fleqn,12pt]{article}
\begin{document}\pagenumbering{arabic}\setcounter{page}{1}\pagestyle{plain}
\baselineskip=16pt

\thispagestyle{empty}
\rightline{MSUMB 97-08, October 1997} 
\vspace{1cm}

\begin{center}
{\Large\bf Two-Parameter Differential Calculus on the $h$-Exterior Plane}
\footnote{This work was supported in part by {\bf T. B. T. A. K.} the 
Turkish Scientific and Technical Research Council. }
\end{center}

\vspace{0.5cm}
\begin{center} Sultan A. Celik 
\footnote{E-mail address: celik@yildiz.edu.tr} \\
Yildiz Technical University, Department of Mathematics, \\
80270 Sisli, Istanbul, TURKEY. 
\end{center}
\begin{center} 
Salih Celik \footnote{E-mail address: 
  scelik@fened.msu..edu.tr}\\
Mimar Sinan University, Department of Mathematics, \\
80690 Besiktas, Istanbul, TURKEY. 
\end{center}

\vspace{2cm}
{\bf Abstract}

We construct a two-parameter covariant differential calculus on the quantum 
$h$-exterior plane. We also give a deformation of the two-dimensional 
fermionic phase space.

\vfill\eject\noindent
{\bf 1. Introduction}

\noindent
A possible approach for quantum groups$^1$ is constructed by 
deforming the coordinates of a space with noncommuting coordinates. 
If we consider the endomorphisms which preserve the algebraic properties 
of the algebra of coordinates then the quantum group structure appears. 

In Ref. 2 Connes showed that derivatives and differentials corresponding 
to these noncommuting coordinates can be defined. He considered the 
differential algebra for noncommutative algebras. It is known, from the 
work of Woronowicz$^3$, that one can define a consistent differential 
calculus on the noncommutative space of a quantum group. 

A one-parameter covariant differential calculus on the $h$-exterior 
plane was introduced in Ref. 4. In this paper a two-parameter 
differential calculus for the noncommutative dual plane is presented. 
It provides a general construction for finding all commutation relations 
among algebra elements, differentials and derivatives. We also give a 
two-parameter deformation of the two-dimensional fermionic phase space and, 
for a special case, we compared the deformed fermionic phase space algebra 
with Ref. 5. We have seen that this deformed algebra, which is obtained 
for a special case, is slightly different from Ref. 5. 

\noindent
{\bf 2. A Differential Calculus on $h$-Exterior Plane} 

\noindent
In this section we formulate a two-parameter differential calculus on the 
quantum $h$-exterior (dual) plane. Let us begin with the basic notations 
concerning the quantum planes. 

The quantum $h'$-plane is defined as an associative algebra whose 
elements $x$, $y$ obey the relation$^{6-8}$ 
$$ xy = yx + h' y^2, \eqno(1)$$
where $h'$ is a complex deformation parameter. The algebra of $h'$-polynomials 
will be called the algebra of functions on the quantum plane and will be 
denoted by ${\cal A}_{h'}$. 

The $h$-exterior plane, denoted by $\Lambda {\cal A}_h$, is an associative 
algebra equipped with the odd (fermionic) generators $\theta$, $\phi$ 
satisfying $^7$ the commutation relations 
$$ \theta^2 = h \theta \phi, \qquad \phi^2 = 0, $$
$$ \theta \phi + \phi \theta = 0 \eqno(2)$$
where $h$ is again a complex deformation parameter. The algebra 
$\Lambda{\cal A}_h$ is a graded algebra with the usual grading induced by 
deg $\theta = 1 = $deg $\phi$. 

In Ref. 9, by formulating a differential calculus on the 
quantum (hyper-) plane, Wess and Zumino showed that the differentials of 
coordinates of the quantum plane can be identified with the coordinates of 
the quantum exterior plane. Since, in the classical limits (i.e. 
$h \longrightarrow 0$, $h' \longrightarrow 0$) $x$ commutes with $y$ 
and $\theta$ anticommutes with $\phi$, as an alternative, 
$x$ and $y$ can be identified with the differentials of $\theta$ and 
$\phi$, respectively, as in Ref. 10. We shall set up a two-parameter 
differential calculus on $\Lambda {\cal A}_h$ interpreting $x$ and $y$ 
as the noncommuting analogues of ${\sf d}\theta$ and ${\sf d}\phi$, 
respectively: 
$$ x = {\sf d}\theta, \quad y = {\sf d}\phi. \eqno(3)$$

To establish a two-parameter differential calculus on the quantum 
$h$-exterior plane, we assume that the commutation relations among the 
coordinates and their differentials are of the form 
$$\Theta^i {\sf d}\Theta^j = C^{ij}{}_{kl} {\sf d}\Theta^k \Theta^l \eqno(4)$$
where the entries $C^{ij}{}_{kl}$ are complex numbers, and 
$\Theta^1 = \theta$, $\Theta^2 = \phi$, etc. 
We would like to describe the entries $C^{ij}{}_{kl}$. 
For this, it is desirable to define an exterior derivative operator {\sf d} 
satisfying the following properties: 

{\bf (i)} {\sf d} is nilpotent, i.e.  
$${\sf d}^2 = 0, \eqno(5)$$

{\bf (ii)} {\sf d} satisfies the graded Leibniz rule, i.e. 
$${\sf d}(f g) = ({\sf d} f) g + (- 1)^{\mbox{deg} f} f ({\sf d} g). \eqno(6)$$

Now from the consistency conditions 
$$ {\sf d}(\theta \phi + \phi \theta) = 0, \qquad 
   {\sf d}(\phi^2) = 0, \qquad {\sf d}(\theta^2 - h \theta \phi) = 0 
\eqno(7)$$
and the derivation of (4) we obtain the matrix 
$$C = 
   \left(\matrix{ 1 & (t - 1) h   & (1 - t) h & (1 - t) h h' \cr 
                  0 &    t        &   1 - t   & (1 - t) h' \cr
                  0 &  1 - t      &     t     & (t - 1) h' \cr 
                  0 &    0        & 0         & 1  \cr }\right) 
  = \hat{R} \eqno(8)$$
where $t$ is a number. The number $t$ can be determined by checking 
the relations among cubic monomials: 
$$ (\theta \phi) {\sf d}\theta = (1 - t) {\sf d}\theta (\phi \theta) + 
   [(t - 1) h' - (t^2 + t - 1) h] {\sf d}\phi (\phi \theta) $$
$$(\phi \theta) {\sf d}\theta = (t^2 + t - 1) {\sf d}\theta (\phi \theta) + 
  (1 - t) (h' - h) {\sf d}\phi(\phi \theta) $$ 
which uniquely constrain $t$ to be equal to 0. The final result is 
$$ \hat{R} 
 = \left(\matrix{ 1 & - h & h &   h h'  \cr 
                  0 &   0 & 1 &   h' \cr
                  0 &   1 & 0 & - h' \cr 
                  0 &   0 & 0 &   1  \cr }\right) \eqno(9)$$
which is given in Ref. 11 (and also Ref. 7). So, we can write the commutation 
relations of coordinates and their differentials as follows: 
$$\Theta^i \Theta^j = - \hat{R}^{ij}{}_{kl} \Theta^k \Theta^l, \eqno(10)$$
$${\sf d}\Theta^i {\sf d}\Theta^j = \hat{R}^{ij}{}_{kl} {\sf d}\Theta^k 
  {\sf d}\Theta^l.  \eqno(11)$$

Now let us denote the partial derivatives with respect to $\theta$ and 
$\phi$, respectively, by 
$$\partial_1 = {\partial\over {\partial \theta}}, \qquad 
  \partial_2 = {\partial\over {\partial \phi}}, \eqno(12)$$
where 
$$\partial_i \Theta^j = \delta^j{}_i. \eqno(13)$$
Assuming the deformed (graded) Leibniz rule for partial derivatives 
$$\partial_i(f g) = (\partial_i f)g + 
  (-1)^{\mbox{deg} f} O^l{}_i(f) \partial_l g   \eqno(14)$$
where 
$$O^l{}_i(\Theta^j) = \hat{R}^{jl}{}_{ik} \Theta^k, \eqno(15)$$
one arrives at 
$$\partial_i \Theta^j = \delta^j{}_i - 
  \hat{R}^{jl}{}_{ik} \Theta^k \partial_l. \eqno(16)$$
Similarly, we have 
$$\partial_i \partial_j = - \hat{R}_{ji}{}^{kl} \partial_l \partial_k, 
  \eqno(17)$$
$$\partial_i {\sf d}\Theta^j = \hat{R}^{jk}{}_{il} {\sf d}\Theta^l \partial_k. 
  \eqno(18)$$

We now order the commutation relations among algebra elements, differentials 
and derivatives as follows: 

\noindent
{\bf (a)} The commutation relations of variables 
$$\theta^2 = h \theta \phi, \quad \phi^2 = 0, $$
$$\theta \phi + \phi \theta = 0. \eqno(19)$$
{\bf (b)} The relations of differentials 
$$x y = y x + h' y^2. \eqno(20)$$
{\bf (c)} The commutation relations of the differentials with the variables 
$$\theta x = x \theta - h (x \phi - y \theta) + h h' y \phi, $$
$$\theta y = y \theta + h' y \phi, $$
$$\phi x = x \phi - h' y \phi, \eqno(21)$$
$$\phi y = y \phi. $$
{\bf (d)} The commutation relations among $\partial$ and $\Theta$ 
$$\partial_\theta \theta = 1 - \theta \partial_\theta + 
  h \phi \partial_\theta, $$
$$\partial_\theta \phi = - \phi \partial_\theta, $$
$$\partial_\phi \theta = - \theta \partial_\phi - h \theta \partial_\theta 
  - h' \phi \partial_\phi - h h' \phi \partial_\theta,   \eqno(22)$$
$$\partial_\phi \phi = 1 - \phi \partial_\phi + 
  h' \phi \partial_\theta. $$
{\bf (e)} The relations of derivatives 
$$\partial_\theta^2 = 0, \quad 
  \partial_\phi^2 = h' \partial_\theta \partial_\phi, $$
$$\partial_\theta \partial_\phi + \partial_\phi \partial_\theta = 0. 
  \eqno(23)$$
{\bf (f)} The commutation relations between derivatives and differentials 
$$\partial_\theta x = x \partial_\theta - h y \partial_\theta, $$
$$\partial_\theta y = y \partial_\theta ,$$
$$\partial_\phi x = x \partial_\phi + h' x \partial_\theta + 
   h y \partial_\phi + h h' y \partial_\theta, $$
$$\partial_\phi y = y \partial_\phi - h' y \partial_\theta. \eqno(24)$$

Note that, it is easy to see that when $h = h'$, this calculus go back to 
those of the one-parameter calculus in Ref. 4. Also, one can shown that 
this two-parameter differential calculus on the $h$-exterior plane (19)-(24) 
is covariant under the action of $GL_{h,h'}(2)$. The quantum group 
$GL_{h,h'}(2)$ was studied in Ref. 7. 

We now combine the relations (19), (22) and (23) and denote the algebra 
generated by the fermionic coordinates $\theta$, $\phi$ and the fermionic 
derivatives $\partial_\theta$, $\partial_\phi$ by ${\cal B}_{h,h'}$.

\noindent
{\bf 3. A Deformation of Fermionic Phase Space}

\noindent
We know that the natural definition of the fermionic momenta is to simply 
identify the fermionic derivatives $\partial_\theta$ and $\partial_\phi$ 
with the fermionic momenta $\pi_\theta$ and $\pi_\phi$. But 
the hermiticity of the fermionic coordinates and the fermionic momenta 
is not compatible with the relations (22). In fact, the only problem is in 
the third relation of (22). In short, the algebra ${\cal B}_{h,h'}$ cannot 
be interpreted as a two-parameter deformation of the fermionic phase space 
algebra. However, if we choose 
$$\overline{h} = - h, \qquad \overline{h'} = - h', \eqno(25)$$
then the algebra ${\cal B}_{h,h'}$ admits the following involution: 
$$\theta^+ = \theta + h \phi, \qquad \phi^+ = \phi \eqno(26)$$
and 
$$\partial_\theta^+ = \partial_\theta, \qquad 
  \partial_\phi^+ = \partial_\phi + h' \partial_\theta. \eqno(27)$$
The above involution allows us to define the hermitean operators 
$$\hat{\theta} = \theta + {h\over 2} \phi, \qquad \hat{\phi} = \phi,\eqno(28)$$
and 
$$\hat{\pi}_\theta = \partial_\theta, \qquad 
  \hat{\pi}_\phi = \partial_\phi + {{h'}\over 2} \partial_\theta. \eqno(29)$$
Then we have 
$$\hat{\theta}^2 = h \hat{\theta} \hat{\phi}, \qquad 
  \hat{\theta} \hat{\phi} + \hat{\phi} \hat{\theta} = 0, \qquad 
  \hat{\phi}^2 = 0, $$
$$\hat{\pi}_\theta^2 = 0, \qquad 
  \hat{\pi}_\theta \hat{\pi}_\phi + \hat{\pi}_\phi \hat{\pi}_\theta = 0, \qquad 
   \hat{\pi}_\phi^2 = h' \hat{\pi}_\theta \hat{\pi}_\phi,  $$
$$\hat{\pi}_\theta \hat{\theta} = 1 - \hat{\theta} \hat{\phi} + 
   h \hat{\phi} \hat{\pi}_\theta, \qquad 
  \hat{\pi}_\theta \hat{\phi} = - \hat{\phi} \hat{\pi}_\theta, \eqno(30)$$
$$\hat{\pi}_\phi \hat{\theta} = - \hat{\theta} \hat{\pi}_\phi - 
   h \hat{\theta} \hat{\pi}_\theta - h' \hat{\phi} \hat{\pi}_\phi + 
   {1\over 2} (h^2 + h'^2) \hat{\phi} \hat{\pi}_\theta +{1\over 2} (h + h'), $$
$$ \hat{\pi}_\phi \hat{\phi} = 1 - \hat{\phi} \hat{\pi}_\phi + 
   h' \hat{\phi} \hat{\pi}_\theta. $$
It is easy to see that when $h = h'$, this algebra coincides with those of 
Ref. 5. One can also directly verify the compatibility of (30) with 
hermiticity of $\hat{\theta}$, $\hat{\phi}$, $\hat{\pi}_\theta$ and 
$\hat{\pi}_\phi$. 

Note that if we take $h' = - h$ then we obtain the algebra 
$$\hat{\theta}^2 = h \hat{\theta} \hat{\phi}, \qquad 
  \hat{\theta} \hat{\phi} + \hat{\phi} \hat{\theta} = 0, \qquad 
  \hat{\phi}^2 = 0, $$
$$\hat{\pi}_\theta^2 = 0, \qquad 
 \hat{\pi}_\theta \hat{\pi}_\phi + \hat{\pi}_\phi \hat{\pi}_\theta = 0, \qquad 
   \hat{\pi}_\phi^2 = h \hat{\pi}_\phi \hat{\pi}_\theta,  $$
$$\hat{\pi}_\theta \hat{\theta} = 1 - \hat{\theta} \hat{\phi} + 
   h \hat{\phi} \hat{\pi}_\theta, \qquad 
  \hat{\pi}_\theta \hat{\phi} = - \hat{\phi} \hat{\pi}_\theta, \eqno(31)$$
$$\hat{\pi}_\phi \hat{\theta} = - \hat{\theta} \hat{\pi}_\phi - 
   h (\hat{\theta} \hat{\pi}_\theta - \hat{\phi} \hat{\pi}_\phi) + 
   h^2 \hat{\phi} \hat{\pi}_\theta, $$
$$ \hat{\pi}_\phi \hat{\phi} = 1 - \hat{\phi} \hat{\pi}_\phi - 
   h \hat{\phi} \hat{\pi}_\theta. $$
This gives a one-parameter deformation of the fermionic phase space algebra 
which may be used to study the two-dimensional quantum fermionic space. We 
also note that the algebra (31) is slightly different from Ref. 5. Thus 
we have also, for one-parameter case, a deformation of the two-dimensional 
fermionic phase space algebra. 

\noindent
{\bf Note Added}

\noindent
We would like to note that the methods used in this letter are quite 
different from the methods of Ref. 12 which contructs a differential 
calculus on the $h$-deformed plane. 

\baselineskip=14pt
\def\refname{References}

\end{document}